\documentclass[seceqn,dvips]{arxbj}
\usepackage{dcolumn}

\aid{0}
\volume{16}
\issue{2}
\pubyear{2010}
\firstpage{585}
\lastpage{603}
\doi{10.3150/09-BEJ210}

\makeatletter

\newcolumntype{d}[1]{D{.}{.}{#1}}
\newtheorem{Th}{Theorem}

\makeatother

\begin{document}
\begin{frontmatter}

\title{A semiparametric efficient estimator in case-control studies}
\runtitle{Semiparametric efficiency in case-control studies}

\begin{aug}
\author{\fnms{Yanyuan} \snm{Ma}\ead[label=e1]{ma@stat.tamu.edu}\corref{}}
\runauthor{Y. Ma}
\address{Department of Statistics,
Texas A\&M University, College Station, TX 77843, USA.\\ \printead{e1}}
\end{aug}

\received{\smonth{4} \syear{2008}}
\revised{\smonth{4} \syear{2009}}

%
\begin{abstract}
We construct a semiparametric estimator in case-control studies where
the gene and the environment are assumed to be independent. A discrete or
continuous parametric distribution of the genes is assumed in the
model. A discrete distribution of the genes can be used to model
the mutation or presence of certain group of genes.
A continuous distribution allows the distribution
of the gene effects to be in a finite-dimensional parametric
family and can hence be used to model the gene expression levels.
We leave the distribution of the environment totally unspecified.
The estimator is derived through calculating the efficiency score
function in a hypothetical setting where a close approximation to
the samples is random. The resulting estimator is proved to be
efficient in the hypothetical situation. The efficiency of the
estimator is further demonstrated to hold in the case-control
setting as well.
\end{abstract}

%
\begin{keyword}
\kwd{case-control study}
\kwd{gene-environment interaction}
\kwd{logistic regression}
\kwd{semiparametric efficiency}
\end{keyword}

\end{frontmatter}

\section{Introduction} \label{sec:intro}

Case-control designs are frequently implemented in clinical studies
where, instead of taking a random sample of a mixed population of both
cases and non-cases, a fixed number of cases and a fixed number of
controls are randomly sampled from the respective populations of cases
and non-cases. Because the resulting samples are no longer random or
independently and identically distributed (i.i.d.), the
classical large-sample asymptotic theories could fail to apply.
In the literature, two main approaches are taken in order to adapt the
large-sample theory
to the case-control setting. The first approach is highlighted in
Breslow \textit{et al.} (\citeyear{BreslowRobinsWellner2000}), where a modified design of the usual
case-control study is proposed. The resulting random sample is
then linked to the true case-control sample through using results from
McNeney (\citeyear{McNeney1998}), where the similarity between random and non-random sample
asymptotic properties is developed by almost establishing the whole
asymptotic theory under non-i.i.d.~samples.
The second approach is somewhat more direct and is
implicitly used by Rabinowitz (\citeyear{Rabinowitz2000}). Instead of treating the
indicator ($D$) of case/control as a random variable, $D$ is assumed
to be
known and all the calculations are performed conditionally on
$D$. Although it does result in the conditional randomness of the case-control
samples, the resulting data is not really identically
distributed. Specifically, two different distributions are involved
and the large-sample theory is still not available.
Strictly speaking, the asymptotic theory for non-i.i.d.~data
rederived in McNeney (\citeyear{McNeney1998}) also needs to be applied in order to treat
such a combination of two sample cases.

In addition to the complexity arising from a case-control design, the
problem considered in this article is also a semiparametric model problem,
whose efficient estimator has not yet been explored even in the
i.i.d.~data situation. Specifically, the problem is as follows.
Suppose that in the general population, the occurrence of a disease ($D=1$)
follows a logistic model $\operatorname{logit}\{\operatorname{Pr}(D=1)\}=m(G,E)$,
where $G$ represents a person's genetic character and $E$ represents the
environmental elements. Further, suppose that $G$ and $E$ are
independent of
each other and that we are interested in the effect of gene, environment
and their interaction on the disease status. Thus,
$m(g,e)=\beta_c+\beta_1g+\beta_2e+\beta_3 ge$. The parametric
form of the distribution of gene $g$ is assumed to be known as
$q(g,\beta_4)$, where $\beta_4$ is an unknown finite-dimensional parameter.
The distribution of the environment, $\eta(e)$, is unspecified.
A special version of this problem is considered in
Chatterjee and Carroll (\citeyear{ChatterjeeCarroll2005}), where $q(g,{\beta_4})$ is assumed to be a discrete
distribution. There, the authors derived a profile maximum likelihood
estimator for $\beta=(\beta_c, \beta_1,\beta_2,\beta_3,\beta_4)^{\mathrm{T}}$ and
showed that it is root-$N$ consistent, where $N$ is the size of the
combined samples. The estimator is later extended to a more general
framework in Spinka \textit{et al.} (\citeyear{SpinkaCarrollChatterjee2005}).
However, it is not investigated whether the
estimator achieves the optimal semiparametric efficiency.

In this paper, we first establish in Section~\ref{sec:casecontrol}
that the classical semiparametric theory of Bickel \textit{et al.}~(\citeyear{BickelEtAl1993}) is applicable
in general case-control studies, without having to rederive the theory in
parallel or having to resort to the results from McNeney (\citeyear{McNeney1998}).
Such first order asymptotic equivalence between case-control sampling
and random sampling is a new result.
We then proceed to compute the semiparametric
efficient score and construct a semiparametric estimator for
$\beta$ in Section~\ref{sec:estimator}.
The computation is carried out in a hypothetical population described
in Section~\ref{sec:casecontrol}. This differs
from the real population from which the cases and controls are drawn.
Hence, the derivation has its own interest and novelty.
In this section, we also prove that
although the estimation of the
nuisance parameter $\eta$ is bypassed in our estimator, the resulting
semiparametric estimator still achieves the optimal efficiency. The
proof and treatment is rather non-standard. Numerical
examples are included in Section~\ref{sec:simulation} to demonstrate
the performance of the proposed
estimator. The performance of the method
in the discrete gene model is close to that of the method in Chatterjee
and Carroll (\citeyear{ChatterjeeCarroll2005})
and we pointed out the possible equivalence between the two methods in
Section~\ref{sec:discuss}.
Some analytical derivations and technical details are included in the
\hyperref[Appendix]{Appendix}.

\section{Case-control data versus i.i.d.~data}\label{sec:casecontrol}

The samples from a case control study are not random because the
disease status is not random. In general, the design randomly samples
$N_1$ individuals from the case population and $N_0$ from the non-case
population.
However, let us consider
a hypothetical population of interest with infinite population size, in
which the
disease to non-disease ratio is fixed at $\pi=N_1/N_0$. Here, the
reason for introducing the notion of hypothetical population is to be
able to use the
classical semiparametric theory for i.i.d.~data, developed in Bickel
\textit{et al.} (\citeyear{BickelEtAl1993}).
If the sample of size $N=N_0+N_1$ from a case-control study happens to
be a random sample from the hypothetical
population of interest, then we have a size-$N$ i.i.d.~random sample and
the usual semiparametric analysis will apply.
The asymptotic results hold when $N\to\infty$ and $\pi$ stays fixed.

Of course, the problem is that
a random sample of size $N$ from the hypothetical population of interest
does not have to have exactly $N_0$ controls and $N_1$ cases, hence we
cannot immediately equate a case-control sample and a random sample
from the hypothetical population. In
general, the number of controls/cases of a random sample from the
hypothetical population will have a binomial
distribution $N_d^r\sim
\operatorname{Binomial}(N, N_d/N)$, $d=0, 1$, which is very close to a normal
distribution when $N$ is large, that is,
$(N_d^r-N_d)/\sqrt{N\pi(1-\pi)} \to\operatorname{Normal}(0, 1)$
in distribution when $N\to\infty$. Here, the superscript $^{r}$
stands for `random.' Furthermore,
the probability of having $|N_d^r-N_d|>N^{2/3}$ goes to zero when
$N\to\infty$. Thus, we could think of the case-control sample as
obtained by randomly picking a size-$N$ sample from the hypothetical
population of
interest, then deleting a random $\mathrm{o}_p(N^{2/3})$ cases (controls) and
adding a random $\mathrm{o}_p(N^{2/3})$ controls (cases). Or, alternatively, we
can think of
the case-control sample as a random sample of size~$N$, but with a
randomly chosen $\mathrm{o}_p(N^{2/3})$ data contaminated in a particular way.
This ``particular'' contamination implies the following three
properties: (i) the contamination
happens only to $\mathrm{o}_p(N)$ of the observations (in the case-control
samples, the contamination in fact only happens to
$\mathrm{o}_p(N^{2/3})$ observations, but, in general, $\mathrm{o}_p(N)$ is already
sufficient for our further analysis);
(ii) the contaminated data is still of
order $\mathrm{O}(1)$, that is, $|X_i^c-X_i|$ is bounded in probability for
$i=1, \dots, N$;
(iii) the zero expectation holds for the contaminated
observations, that is, if an estimating equation for $\beta$ of the
form $\sum_{i=1}^N
f(X_i;\beta)=0$ satisfies $E\{f(X_i;\beta_0)\}=0$, then $E\{
f(X_i^c;\beta_0)\}=0$ as
well. Here, $X_i, i=1, \dots, N$, are i.i.d.~random samples,
the superscript $^c$ stands for `contaminated' and the
subscript $_0$ represents the true parameter value.

When the case-control sample is viewed as a contaminated random sample
from the hypothetical population of interest, the first two
``particular'' properties certainly hold. For the estimator we will
construct, we shall demonstrate that the third property also holds.
Thus, if we can show that the same first order asymptotics
apply to both the i.i.d.~sample of size $N$ and its
contaminated version as long as the three
properties hold, then we can treat the
case-control sample as an i.i.d.~sample.

The argument is as follows.
Assume that we mistakenly treated the contaminated data as i.i.d.~and
obtained an efficient
estimator:
%
\begin{eqnarray}\label{eq:score}
\sum_{i=1}^N S_{\mathrm{eff}}(X_i^c;\beta)=0.
\end{eqnarray}
Here, $S_{\mathrm{eff}}$ is the efficient score function and
its derivation is model-dependent. One obvious
aspect of $S_{\mathrm{eff}}$ worth emphasizing is that the
construction of $S_{\mathrm{eff}}$ does not depend on the observations.
Regardless of the method of
derivation, the efficient score
function $S_{\mathrm{eff}}$ has the property
$E\{S_{\mathrm{eff}}(X_i;\beta_0)\}=0$. If we had\vspace*{1pt}
the uncontaminated data,
our subsequent estimator for $\beta$ would have been
$\sum_{i=1}^N S_{\mathrm{eff}}(X_i;\beta)=0$. Working with the
contaminated
data, (\ref{eq:score}) is the estimating equation we really have.
Suppose that $\hat\beta$ solves (\ref{eq:score}). We then have
\begin{eqnarray*}
0=\sum_{i=1}^N S_{\mathrm{eff}}(X_i^c;\hat\beta)
=\sum_{i=1}^N S_{\mathrm{eff}}(X_i^c;\beta_0)+\sum_{i=1}^N
\frac{\partial
S_{\mathrm{eff}}(X_i^c;\beta^*)}{\partial\beta^{\mathrm{T}}}(\hat\beta-\beta_0),
\end{eqnarray*}
therefore,
%
\begin{eqnarray}\label{eq:expand}
-N^{-1}\Biggl\{\sum_{i=1}^N \frac{\partial
S_{\mathrm{eff}}(X_i^c;\beta^*)}{\partial\beta^{\mathrm{T}}}\Biggr\}
\sqrt{N}(\hat\beta-\beta_0)=N^{-1/2}\sum_{i=1}^N S_{\mathrm{eff}}(X_i^c;\beta_0),
\end{eqnarray}
where $\beta^*$ lies on the line connecting $\beta_0$ and $\hat\beta$.
Note that in our ``particular'' contamination requirement,
only $\mathrm{o}_p(N)$ terms yield a different $X_i$ from $X_i^c$ (requirement (i))
and, for each $X_i^c\ne X_i$, the difference is $\mathrm{O}_p(1)$ (requirement (ii)),
so we have
%
\begin{eqnarray}\label{eq:part1}
N^{-1}\Biggl\{\sum_{i=1}^N \frac{\partial
S_{\mathrm{eff}}(X_i^c;\beta^*)}{\partial\beta^{\mathrm{T}}}\Biggr\}
&=&N^{-1}\Biggl\{\sum_{i=1}^N \frac{\partial
S_{\mathrm{eff}}(X_i;\beta^*)}{\partial\beta^{\mathrm{T}}}\Biggr\}+\mathrm{o}_p(1)\nonumber
\\[-8pt]\\[-8pt]
&=&E\Biggl\{\frac{\partial
S_{\mathrm{eff}}(X_i;\beta_0)}{\partial\beta^{\mathrm{T}}}\Biggr\}+\mathrm{o}_p(1).\nonumber
\end{eqnarray}
From the third ``particular'' property, we have
$E\{S_{\mathrm{eff}}(X_i^c;\beta_0)\}=0$
(we will prove that this property holds for the case-control data in Section~\ref{sec:estimator}). In conjunction with
the fact that only $\mathrm{o}_p(N)$ of the terms $S_{\mathrm{eff}}(X_i^c;\beta_0)
-S_{\mathrm{eff}}(X_i;\beta_0) $ are non-zero,
we can further obtain
%
\begin{eqnarray}\label{eq:part2}
N^{-1/2}\sum_{i=1}^N S_{\mathrm{eff}}(X_i^c;\beta_0)
=N^{-1/2}\sum_{i=1}^N S_{\mathrm{eff}}(X_i;\beta_0) +\mathrm{o}_p(1).
\end{eqnarray}
The detailed argument of (\ref{eq:part2}) is the following. Suppose
for the first $l=\mathrm{o}_p(N)$
observations, $X_i^c\ne X_i$. Then we have
\begin{eqnarray*}
&&N^{-1/2}\sum_{i=1}^N S_{\mathrm{eff}}(X_i^c;\beta_0)
\\
&&\quad=N^{-1/2}\sum_{i=1}^N S_{\mathrm{eff}}(X_i;\beta
_0)+N^{-1/2}\sum_{i=1}^l
\{S_{\mathrm{eff}}(X_i^c;\beta_0)-S_{\mathrm{eff}}(X_i;\beta_0)\}\\
&&\quad=N^{-1/2}\sum_{i=1}^N S_{\mathrm{eff}}(X_i;\beta
_0)+(N/l)^{-1/2}l^{-1/2}\sum_{i=1}^l
\{S_{\mathrm{eff}}(X_i^c;\beta_0)-S_{\mathrm{eff}}(X_i;\beta_0)\}.
\end{eqnarray*}
Note that $S_{\mathrm{eff}}(X_i^c;\beta_0)-S_{\mathrm{eff}}(X_i;\beta_0)$ has mean
zero, hence $l^{-1/2}\sum_{i=1}^l\{S_{\mathrm{eff}}(X_i^c;\beta_0)-S_{\mathrm{eff}}(X_i;\beta_0)\}=\mathrm{O}_p(1)$. From
$l=\mathrm{o}_p(N)$, we obtain the result in (\ref{eq:part2}) immediately.
Thus, plugging (\ref{eq:part1}) and (\ref{eq:part2}) into
(\ref{eq:expand}),
we obtain
\begin{eqnarray*}
-E\biggl\{\frac{\partial
S_{\mathrm{eff}}(X_i;\beta_0)}{\partial\beta^{\mathrm{T}}}\biggr\}
\sqrt{N}(\hat\beta-\beta_0)=N^{-1/2}\sum_{i=1}^N S_{\mathrm{eff}
}(X_i;\beta_0)+\mathrm{o}_p(1).
\end{eqnarray*}
The above display is exactly the first order asymptotic expansion of
the estimator for $\beta$ if we had performed the estimation procedure
on the uncontaminated data. Thus, we have demonstrated that the
estimator obtained from contaminated data performs as well as the one
obtained from uncontaminated data in terms of first order asymptotic
properties.
Note that the efficient estimator can be replaced by a consistent
estimator, say, a general $S$ instead of $S_{\mathrm{eff}}$, as
long as
$E(S|D=d)=0$ holds for $d=0, 1$. This ensures that
$E\{S(X_i^c)\}=0$ as long as $E\{S(X_i)\}=0$
(shown in Section~\ref{sec:estimator}),
so the above derivation will still carry through.
Hence, the asymptotic property of the estimator using the
contaminated data is indeed the same as if we had the uncontaminated
data. Thus,
the case-control data can be treated as i.i.d.~data and
we can achieve the same efficiency as when the data was indeed
i.i.d. In other words, a semiparametric estimator using contaminated
data is
at least as efficient as one using the uncontaminated data.

One question still remains: can we do even better than in the i.i.d.~data
case?
In fact, since case-control sampling is designed to be an efficient
way to collect covariate information, it seems to contain more
information than a random sample. However, we claim that
for asymptotically linear estimators of the form
\begin{eqnarray*}
\sqrt{N}(\hat\beta-\beta_0)=\frac{1}{\sqrt{N}}\sum_{i=1}^N
\psi(X_i^c;\beta_0)+\mathrm{o}_p(1),
\end{eqnarray*}
where $E\{\psi(X_i^c;\beta_0)|d\}
=0$, the efficiency in parameter estimation cannot be
further improved by taking into account the special sampling
procedure. This is
because otherwise, we could have obtained a better
estimator for the i.i.d.~sample as well, by replacing $X_i^c$ with
$X_i$.
The detailed derivation is the same as in the
above paragraph, where the condition $E\{\psi(X_i^c;\beta_0)|d\}
=0$ implies $E\{\psi(X_i;\beta_0)|d\}=0$ for
case-control data, which
ensures~\mbox{$E\{\psi(X_i^c;\beta_0)\}=E\{\psi(X_i;\beta_0)\}=0$}.
Of course, if the condition $E(\psi|d)=0$ is not satisfied, the
argument does
not work. However, we now show that if $\psi$ achieves the optimal variance
for the case-control data $X_i^c$, then it has to satisfy
$E\{\psi(X_i^c;\beta_0)|d\}=0$.

First,\vspace*{1pt} $E\{\partial E(\psi|D)/\partial\beta\}
=\partial E(\psi)/\partial\beta=0$ because the probability density
function (p.d.f.) of $D$ does not contain $\beta$.
If we let
$\tilde\psi(X_i^c)=\psi(X_i^c)-E\{\psi(X_i^c)|d\}$, then
$E\{\tilde\psi(X_i^c)\}=0$ and
$E\{\partial\tilde\psi(X_i^c)/\partial\beta\}=E\{\partial\psi
(X_i^c)/\partial\beta\}$.
If $E\{\psi(X_i^c)|d\}\ne0$, then we can obtain
\begin{eqnarray*}
\operatorname{var}\{\psi(X_i^c)\}
&=&E[\operatorname{var}\{\psi
(X_i^c)|D\}]+\operatorname{var}[E\{\psi(X_i^c)|D\}]
=\operatorname{var}\{\tilde\psi(X_i^c)\}+\operatorname{var}[E\{\psi(X_i^c)|D\}]
\\
&>&\operatorname{var}\{\tilde\psi(X_i^c)\},
\end{eqnarray*}
which, together with
$E\{\partial\tilde\psi(X_i^c)/\partial\beta\}=E\{\partial\psi
(X_i^c)/\partial\beta\}$,
contradicts the fact that $\psi(X_i^c)$ is optimal.

In summary, we have shown that the case control samples can be treated
as if they were i.i.d.~and all the first order asymptotic results for
i.i.d.~data will be inherited for case-control data as well. We can see that
the above establishment is similar to the development in Breslow \textit{et al.}~(\citeyear{BreslowRobinsWellner2000}).
However, one prominent difference is that in
Breslow \textit{et al.}~(\citeyear{BreslowRobinsWellner2000}), the case-control sample is viewed as the result of a biased
sampling procedure with fixed subsample size, hence they cannot use
the classical semiparametric theory for i.i.d.~data, but have to refer
to McNeney (\citeyear{McNeney1998}) for the theoretical properties, where the whole
semiparametric theory for fixed-size subsamples is established in
parallel to the
i.i.d.~framework. Here, through introducing the notion of
hypothetical population and by analyzing the first order equivalence
between a
random sample and a sample with fixed-size subsamples, we can easily
contain the case-control problem in the usual i.i.d.~model framework.
The derivation is much simpler and more elegant.
Thus, in the remainder of the paper, we ignore the case-control nature
of the data and
proceed with our analysis by pretending the data is i.i.d.~from the
aforementioned hypothetical population of interest.

\section{A semiparametric efficient estimator}\label{sec:estimator}
\subsection{Geometric approach}\label{sec:geometry}
A random sample from the hypothetical
population of interest has p.d.f.
%
\begin{eqnarray}\label{eq:pdf}
p(g,e,d; \beta,\eta)
&=&p_D(d)p_{G,E|D}(g,e|d)=p_D(d)p_{G,E|D}^t(g,e|d)\nonumber
\\
&=&p_D(d)p_{G}^t(g)p_{E}^t(e)p^t_{D|G,E}(d|g,e)/p_{D}^t(d)
\\
&=&\frac{N_d}{N}\frac{q(g)\eta(e)H(d,g,e)}{p_D^t(d)}.\nonumber
\end{eqnarray}
Here, the superscript $^t$ stands for the p.d.f.~in the
true population, whereas
expressions without
superscripts, including various p.d.f.'s and expectation $E$,
are quantities in the hypothetical population of
interest;
$\eta(e)=p_{E}^t(e)$ is the unknown infinite-dimensional parameter
and
\begin{eqnarray*}
H(d,g,e;\beta)&=&\exp[d\{m(g,e)\}]/[1+\exp\{m(g,e)\}]
\\
&=&\exp\{d(\beta_c+\beta_1g+\beta_2e+\beta_3ge)\}/\{1+\exp(\beta_c+\beta_1g+\beta_2e+\beta_3ge)\};
\\
p_D^t(d; \beta,\eta)
&=&\int q(g,\beta_4)\eta(e)H(d,g,e;\beta)\,\mathrm{d}\mu(g)\,\mathrm{d}\mu(e).
\end{eqnarray*}
We recognize that estimating the finite-dimensional parameter $\beta$
in the presence of
an infinite-dimensional nuisance parameter $\eta$, using an
i.i.d.~sample of size $N=N_0+N_1$ from a
hypothetical population of interest, with the p.d.f.~of a random observation
given by (\ref{eq:pdf}), is a classical semiparametric
problem. Therefore, we implement the semiparametric estimation methods
to derive the semiparametric efficient estimator. The approach we take
is geometric, first introduced in Bickel \textit{et al.}~(\citeyear{BickelEtAl1993}).
Because the
general approach and related concepts have been nicely described in several
recent papers including Tsiatis and Ma (\citeyear{TsiatisMa2004}),
Allen \textit{et al.}~(\citeyear{AllenSattenTsiatis2005}),
Ma \textit{et al.}~(\citeyear{MaGentoTsiatis2005}) and
Ma and Tsiatis (\citeyear{MaTsiatis2006}), here, we only briefly outline
the general approach and the definition of the relevant concepts,
referring the reader to these papers for more detailed descriptions.

In general semiparametric problems,
one approach to construct estimators for $\beta$ is
to obtain some influence function $\phi(X_i;\beta,\eta)$
which is subsequently used to form estimating equations for $\beta$ in
the form of $\sum_{i=1}^N \phi(X_i;\beta,\eta)=0$. Here,
$X_i=(G_i,E_i,D_i),~i=1,\dots,N$, are i.i.d.~observations. The
solution of the
estimating equation, $\hat\beta$, is subsequently a semiparametric
estimator and its variance has been established to be
equal to the variance of $\phi(X_i;\beta,\eta)$. Consequently, the
optimal estimator among the
class of all such estimators is the one whose
influence function has the smallest variance. This is usually referred
to as the \textit{semiparametric efficient estimator}.

The geometric approach considers the space in which all influence
functions belong. Specifically, one
considers a Hilbert space $\mathcal{H}$ which consists of all
zero-mean measurable functions
with finite variance and the same dimension as $\beta$.
The inner product in $\mathcal{H}$ is defined as the covariance.
The Hilbert space $\mathcal{H}$ is further decomposed into two spaces, the
nuisance tangent space $\Lambda$ and its orthogonal complement
$\Lambda^\perp$.

To understand the nuisance tangent space $\Lambda$, consider first
the case where the nuisance parameter, denoted $\gamma$, is finite-dimensional.
Then, the nuisance score function, $S_\gamma=\partial\log
p(X_i;\beta,\gamma)/\partial\gamma$, spans a linear space, which is denoted
$\Lambda$. In the case of the infinite-dimensional nuisance parameter
$\eta$, the
corresponding $\Lambda$ is defined as the mean squared closure of the
span of all the nuisance score functions $S_\gamma$, where
$p(X_i;\beta,\gamma)$ is any parametric submodel of
$p(X_i;\beta,\eta)$. The orthogonal complement of $\Lambda$ in $\mathcal{H}$ is subsequently defined as~$\Lambda^\perp$.

Any function in $\Lambda^\perp$ can be properly
normalized to obtain a valid influence function. On the other hand,
every influence function is a function in $\Lambda^\perp$. Among all these
functions, the
projection of the score function $S_\beta=\partial\log
p(X_i;\beta,\gamma)/\partial\beta$ results in the efficient influence
function. If we denote the projection by $S_{\mathrm{eff}}$,
then the corresponding
optimal variance is $\operatorname{var}(S_{\mathrm{eff}})^{-1}$. The projection $S_{\mathrm{eff}}$
is usually called the \textit{efficient score function}.

Hence, the geometric approach converts the problem of searching for
efficient semiparametric estimators to the problem of calculating
$S_{\mathrm{eff}}$.

\subsection{Construction of the estimator}\label{sec:construct}

Following the description in Section~\ref{sec:geometry}, we obtain
the efficient score function $S_{\mathrm{eff}}$.
Viewing the sample as random from the hypothetical
population, the p.d.f.~in (\ref{eq:pdf}) is no longer in a simple
multiplicative form,
in that the nuisance parameter appears both in the numerator and in
the integral in the denominator.
Since this implies that the nuisance tangent space is not
automatically orthogonal to the score functions,
the related computation for the nuisance tangent space
and associated objects
is more involved. In addition, one needs to be aware that the
calculation should be carried out with respect to the hypothetical
population, hence quantities such as $p_G^t, p_E^t, p_D^t$ need to
be treated with extra care and not confused with
$p_G, p_E, p_D$.
The main steps of the derivation are as follows.
We first calculate the score function $S_\beta$ by taking the derivative
of $\log p(g,e,d;\beta,\eta)$ with respect to $\beta$. This results
in $S_\beta=S-E(S|d)$, where
\begin{eqnarray*}
S=\biggl\{(m'_{\beta_c}\
m'_{\beta_1}\
m'_{\beta_2}\
m'_{\beta_3})
\biggl(d-1+\frac{1}{1+\mathrm{e}^m}\biggr)\quad
\frac{q'_{\beta_4}(g,{\beta_4})^{\mathrm{T}}}{q(g,{\beta_4})}\biggr\}^{\mathrm{T}}.
\end{eqnarray*}
We then calculate the two spaces $\Lambda, \Lambda^\perp$ by replacing
$\eta$ in (\ref{eq:pdf}) with a finite-dimensional parameter $\gamma$,
taking the derivative of $\log p(g,e,d;\beta,\gamma)$ with respect to
$\gamma$ to obtain $S_\gamma$, hypothesizing a space of all such
$S_\gamma$ and proving that $\Lambda$ is equivalent to this space.
The results are
\begin{eqnarray*}
\Lambda&=&[h(e)-E\{h(e)|d\}{}\dvt{} \forall h(e) \mbox{ such that } E^t(h)=0]
=[h(e)-E\{h(e)|d\}{}\dvt{} \forall h(e)],\\
\Lambda^\perp&=&[h(g,e,d){}\dvt{} E(h|e)=E\{E(h|d)|e\}].
\end{eqnarray*}
We finally project the score vector $S_\beta$ onto $\Lambda^\perp$ to
obtain $S_{\mathrm{eff}}=S_\beta
-f(e)+E(f|d)=S-E(S|d)-f(e)+E(f|d)$, where
$f(e)-E(f|d)$ represents the projection of $S_\beta$ onto $\Lambda$.
The details of the derivation can be found in the \hyperref[Appendix]{Appendix}.
Note that this form of $S_{\mathrm{eff}}$ implies that $E\{
S_{\mathrm{eff}}(X)|d\}=0$.
When $X$ is replaced by $X^c$, the non-random case-control sample, we
still have $E\{S_{\mathrm{eff}}(X^c)|d\}=0$ because the design
itself guarantees
that the only non-random component is~$d$, which is held constant.
Thus, viewing $X^c$ as a special contaminated version of $X$,
we still have $E\{S_{\mathrm{eff}}(X^c)\}=0$, which is required
in Section~\ref{sec:casecontrol}.

From the \hyperref[Appendix]{Appendix}, we can further write
%
\begin{eqnarray}\label{eq:effscore}
S_{\mathrm{eff}}
&=&S-E(S|e)+(-1)^d\{a(0)-a(1)\}w(e,1-d),
\end{eqnarray}
where $a(0)-a(1)=E(f|D=0)-E(S|D=0)-E(f|D=1)+E(S|D=1)$.

In terms of the calculation of $S_{\mathrm{eff}}$,
note that $S$, $E(S|e)$ and $w$, as given in (\ref{eq:w}),
are all functions with parameters $\beta$ and $p_D^t(d)$ only. Hence,
as long as we can calculate $p_D^t(d)$,
we will have the ability to evaluate $S$, $E(S|e)$ and $w$.
The computation of $a(0)-a(1)$ requires further arguments.

In the following, we first obtain an approximation of $p_D^t(d)$, then
pursue the estimation of $a(0)-a(1)$.
To estimate $p_D^t(d)$, using $p_E(e)$ to denote the probability density
function of $e$ in the hypothetical population, we observe that
\begin{eqnarray*}
N_d&=&Np_D(d)=\int Np_{D,E}(d,e) \,\mathrm{d}\mu(e)
=\int Np_E(e)p_{D,G|E}(d,g|e)\,\mathrm{d}\mu(g)\,\mathrm{d}\mu(e)\\
&=&\int
Np_E(e)\frac{ \int N_dq(g,{\beta_4})H(d,g,e)\,\mathrm{d}\mu(g)/p_D^t(d)}
{ \sum_d\int N_dq(g,{\beta_4})H(d,g,e)\,\mathrm{d}\mu(g)/p_D^t(d)}\,\mathrm{d}\mu(e)\\
&=&E_e\biggl\{\frac{N\int N_dq(g,{\beta_4})H(d,g,e)\,\mathrm{d}\mu(g)/p_D^t(d)}
{\sum_d\int N_dq(g,{\beta_4})H(d,g,e)\,\mathrm{d}\mu(g)/p_D^t(d)}\biggr\}.
\end{eqnarray*}
Replacing the moment $E_e$ with its sample moment through averaging
across different observed $e_i$'s, we obtain
%
\begin{eqnarray}\label{eq:pd}
N_d\approx\sum_{i=1}^N \frac{\int N_dq(g,{\beta_4})H(d,g,e_i)\,\mathrm{d}\mu(g)/p_D^t(d)}
{\sum_d\int N_dq(g,{\beta_4})H(d,g,e_i)\,\mathrm{d}\mu(g)/p_D^t(d)}\qquad \mbox{for } d=0,1.
\end{eqnarray}
Note that the above two equations are not independent -- one determines
the other. But, in combination with $p_D^t(0)+p_D^t(1)=1$, we can
estimate $p_D^t(d)$ completely.
Because the only approximation involved in estimating $p_D^t(d)$ is
replacing the mean with a sample mean, the calculation
will produce a root-$N$-consistent
estimator for $p_D^t(0)$ and $p_D^t(1)$. We denote the estimators by
$\hat p_D^t(0)$ and $\hat p_D^t(1)$.
In calculating $N_d$, we write $p(g,e,d)$ as $p_E(e) p_{D,G|E}(d,g|e)$,
instead of directly using the form in (\ref{eq:pdf}). Since $p_E(e)$ is
the p.d.f.~of the environment variable in the hypothetical population,
this enables us to replace the expectation $E_e$ with the average of the samples.\looseness=-1

The estimation of $a(0)-a(1)$ is much more tedious, and involves
an almost brute force calculation of $E(S|d)$ and $E(f|d)$.
If we let $b_0=E(S|D=0), b_1=E(S|D=1), c_0=E(f|D=0)$ and
$c_1=E(f|D=1)$, then
$a(0)-a(1)=b_1-b_0+c_0-c_1$. The calculation of $b_0$ and $b_1$
follows from
\begin{eqnarray*}
b_d
&=&\frac{\int Sp_{D,G,E}(d,g,e)\,\mathrm{d}\mu(g)\,\mathrm{d}\mu(e)}
{\int p_{D,G,E}(d,g,e)\,\mathrm{d}\mu(g)\,\mathrm{d}\mu(e)}
=\frac{\int Sp_E(e)p_{D,G|E}(d,g|e)\,\mathrm{d}\mu(g)\,\mathrm{d}\mu(e)}
{\int p_E(e)p_{D,G|E}(d,g|e)\,\mathrm{d}\mu(g)\,\mathrm{d}\mu(e)}\\
&=&\int p_E(e)\frac{\int SN_dq(g)H(d,g,e)\,\mathrm{d}\mu(g)/p_D^t(d)}
{\sum_d\int N_dq(g)H(d,g,e) \,\mathrm{d}\mu(g)/p_D^t(d)}\,\mathrm{d}\mu(e)
\\
&&{}\Big/
\int p_E(e)\frac{\int N_dq(g)H(d,g,e)\,\mathrm{d}\mu(g)/p_D^t(d)}
{\sum_d\int N_dq(g)H(d,g,e) \,\mathrm{d}\mu(g)/p_D^t(d)}\,\mathrm{d}\mu(e).
\end{eqnarray*}
Since $S$ can be calculated directly, we simply obtain the
approximation of $b_d,d=0,1$, by replacing the mean with sample
mean and plugging in the estimated $p_D^t(d)$:
%
\begin{eqnarray}\label{eq:b0}
\hat b_0&=&\sum_{i=1}^N
\frac{\int S(0,g,e_i)q(g)H(0,g,e_i)\,\mathrm{d}\mu(g)}
{\sum_d\int N_dq(g)H(d,g,e_i) \,\mathrm{d}\mu(g)/\hat p_D^t(d)}\nonumber
\\[-8pt]\\[-8pt]
&&{}\Big/
\sum_{i=1}^N \frac{\int q(g)H(0,g,e_i)\,\mathrm{d}\mu(g)}
{\sum_d\int N_dq(g)H(d,g,e_i) \,\mathrm{d}\mu(g)/\hat p_D^t(d)},\nonumber
\\
\hat b_1&=&\sum_{i=1}^N\frac{\int S(1,g,e_i)q(g)H(1,g,e_i)\,\mathrm{d}\mu(g)}
{\sum_d\int N_dq(g)H(d,g,e) \,\mathrm{d}\mu(g)/\hat p_D^t(d)}\label{eq:b1}\nonumber
\\[-8pt]\\[-8pt]
&&{}\Big/\sum_{i=1}^N\frac{\int q(g)H(1,g,e_i)\,\mathrm{d}\mu(g)}
{\sum_d\int N_dq(g)H(d,g,e) \,\mathrm{d}\mu(g)/\hat p_D^t(d)}.\nonumber
\end{eqnarray}
The calculations of $c_0$ and $c_1$ are a bit more tricky. Since
\begin{eqnarray*}
&&f=E(S|e)+(c_0-b_0)w(e,0)+(c_1-b_1)\{1-w(e,0)\},
\end{eqnarray*}
taking expectation conditional on, say $D=0$, we have
\begin{eqnarray*}
c_0&=&E\{E(S|e)|D=0\}+(c_0-b_0)E\{w(e,0)|D=0\}
\\
&&{}+(c_1-b_1)[1-E\{w(e,0)|D=0\}]
\end{eqnarray*}
or, equivalently, we obtain
\begin{eqnarray*}
c_0-c_1=\frac{E\{E(S|e)|D=0\}-b_0E\{w(e,0)|D=0\}-b_1[1-E\{w(e,0)|D=0\}
]}{1-E\{w(e,0)|D=0\}}.
\end{eqnarray*}
Hence, replacing mean by sample mean and using $\hat p_D^t(d)$,
$c_0-c_1$ is estimated by
%
\begin{eqnarray}\label{eq:c}
\hat c_0-\hat c_1=
\frac{\hat E\{E(S|e)|D=0\}-\hat b_0\hat E\{w(e,0)|D=0\}-\hat
b_1[1-\hat E\{w(e,0)|D=0\}]}{1-\hat E\{w(e,0)|D=0\}},
\end{eqnarray}
where
%
\begin{eqnarray}\label{eq:w0}
\hat E\{w(e,0)|D=0\}\nonumber
&=&\sum_{i=1}^N\frac{w(e_i,0)\int
q(g)H(0,g,e_i)\,\mathrm{d}\mu(g)}{\sum_d\int N_dq(g)H(d,g,e_i) \,\mathrm{d}\mu(g)/\hat p_D^t(d)}\nonumber
\\[-8pt]\\[-8pt]
&&{}\Big/
\sum_{i=1}^N \frac{\int q(g)H(0,g,e_i)\,\mathrm{d}\mu(g)}{\sum_d\int N_dq(g)H(d,g,e_i) \,\mathrm{d}\mu(g)/\hat p_D^t(d)}\nonumber
\end{eqnarray}
and
\begin{eqnarray}\label{eq:w1}
\hat E\{E(S|e)|D=0\}\nonumber
&=&\sum_{i=1}^N\frac{E(S|e_i)\int
q(g)H(0,g,e_i)\,\mathrm{d}\mu(g)}{\sum_d\int N_dq(g)H(d,g,e_i) \,\mathrm{d}\mu(g)/\hat p_D^t(d)}\nonumber
\\[-8pt]\\[-8pt]
&&{}\Big/
\sum_{i=1}^N \frac{\int q(g)H(0,g,e_i)\,\mathrm{d}\mu(g)}{\sum_d\int N_dq(g)H(d,g,e_i) \,\mathrm{d}\mu(g)/\hat p_D^t(d)}.\nonumber
\end{eqnarray}
Similarly to the estimation of $p_D^t(d)$, the only approximation
involved in obtaining $b(0), b(1)$ and $c(0)-c(1)$ is replacing mean by
sample mean, so
$a(0)-a(1)$ is estimated using $\hat a(0)-\hat a(1)=\hat b_1-\hat
b_0+\hat c_0-\hat c_1$ at the root-$N$ rate.

We would like to emphasize that in all of the above calculations,
when we replace the expectation with the sample average, we use the
result that the case-control sample can be treated as a random sample\vspace*{1pt}
from the hypothetical population. Hence, for any function $u(e)$,
the approximation $N^{-1}\sum_{i=1}^Nu(e_i)$
can only be used to replace $\int u(e)p_E(e)\,\mathrm{d}\mu(e)$, not $\int
u(e)\eta(e) \,\mathrm{d}\mu(e)$.

We omitted the parameter $\beta$ in all of the above expressions, in
fact,
$p_D^t(0), p_D^t(1), a(0)-a(1)$ are all functions of $\beta$.
However, if we replace $\beta$ with $\tilde\beta$, an initial
estimator of $\beta$, we will still obtain
$\hat p_D^t(d;\tilde\beta), \hat a(0;\tilde\beta)-\hat
a(1;\tilde\beta)$ that are root-$N$-consistent, as long as
$\tilde\beta-\beta=\mathrm{O}_p(N^{-1/2})$.
The final estimating equation of $\beta$ has the form
%
\begin{eqnarray}\label{eq:esteq}
\sum_{i=1}^N\hat S_{\mathrm{eff}}(x_i;\beta)=\sum
_{i=1}^NS_{\mathrm{eff}}\{
x_i;\beta,\hat p_D^t(d;\tilde\beta), \hat a(0;\hat
p_D^t,\tilde\beta)-\hat a(1;\hat p_D^t,\tilde\beta)\}=0,
\end{eqnarray}
where $x_i$ denotes the $i$th observation $(d_i,g_i,e_i)$.

To summarize the description of the estimator, we outline the
algorithm here:
\begin{enumerate}[Step 5.]
\item[Step 1.]
Pick a starting value $\tilde\beta$ that is root-$N$ consistent.
\item[Step 2.] Solve for $\hat p_D^t(0)$ and $\hat p_D^t(1)=1-\hat p_D^t(0)$ from (\ref{eq:pd}).
\item[Step 3.]
Obtain $\hat b_0$ and $\hat b_1$ from (\ref{eq:b0}) and
(\ref{eq:b1}).
\item[Step 4.]
Obtain $\hat c_0-\hat c_1$ from (\ref{eq:c}) and (\ref{eq:w0}), (\ref{eq:w1}).
\item[Step 5.]
Calculate $S_{\mathrm{eff}}$ using (\ref{eq:effscore}) and
obtain $\hat\beta
$ from
solving (\ref{eq:esteq}).
\end{enumerate}

It is worth pointing out
that in order to carry out Step 1,
we have used a vital assumption that a root-$N$ starting value
$\tilde\beta$ exists. Fortunately, the existence of $\tilde\beta$ is
equivalent to the identifiability of $\beta$ and is already well
established in Chatterjee and Carroll (\citeyear{ChatterjeeCarroll2005}). The starting value
used there, or in Spinka \textit{et al.}~(\citeyear{SpinkaCarrollChatterjee2005}),
can be used to obtain the initial estimator $\tilde\beta$.
Our algorithm here does not require an iteration of
Steps 2--5 upon each update of $\beta$. However, in practice, a more
accurate $\tilde\beta$ can improve the final estimation $\hat\beta$
significantly, hence iterations are almost always implemented.

\subsection{Semiparametric efficiency}\label{sec:efficiency}
If we could use the exact $p_D^t(d;\beta)$ and $a(0;\beta)-a(1;\beta
)$ in
(\ref{eq:esteq}),
then, according to Section~\ref{sec:geometry}, the resulting
estimator for $\beta$ would be an efficient estimator,
with estimation variance
$V=E(S_{\mathrm{eff}}S_{\mathrm{eff}}^{\mathrm{T}})^{-1}$.
To first order, $V$ can be approximated
using $N\{\sum_{i=1}^N\hat S_{\mathrm{eff}}(x_i;\hat\beta
)\hat
S_{\mathrm{eff}}^{\mathrm{T}}(x_i;\hat\beta)\}^{-1}$,
where $\hat\beta$ solves (\ref{eq:esteq}).

We claim that using the estimated $\hat S_{\mathrm{eff}}$ as in
(\ref
{eq:esteq}), we
obtain an estimating equation that yields the same estimator for
$\beta$ as using $S_{\mathrm{eff}}$, in terms of its first
order asymptotic
properties.
\begin{Th}\label{th:thm}
The algorithm in Section~\ref{sec:construct} yields a semiparametric
efficient estimator for $\beta$. That is,
\begin{eqnarray*}
\sqrt{N}(\hat\beta-\beta_0)\to\operatorname{Normal}\{0,\operatorname{var}(S_{\mathrm{eff}})^{-1}\}
\end{eqnarray*}
in distribution when $N\to\infty$ and $N_1/N_0$ is fixed.
\end{Th}

The proof of the theorem contains two main steps. In the first step,
we show the semiparametric efficiency of the estimator if the
observations had been i.i.d. In the second step, we proceed to show
the efficiency in the case-control study using results in Section~\ref{sec:casecontrol}. Rather complex algebra needs to be employed in
the first
step. The proof also involves a split of the data in the final
estimation of $\beta$, and in estimating $p_D^t(d)$ and
$a(0)-a(1)$, mainly for technical convenience.
The details of the proof appear in the \hyperref[Appendix]{Appendix}.

\section{Numerical examples}\label{sec:simulation}
We conducted a small simulation study to demonstrate the performance of the
estimator.
In the first experiment, we generated 500 cases and 500 controls,
where the true environment element $E$ is
$\min(10,X)$ and $X$ is generated from a
log-normal distribution with mean 0 and variance 1. A dichotomous model
of the gene is used, where $G=1$ with probability $\beta_4$ and $G=0$
with probability $1-\beta_4$. This kind of model for $q(g,\beta_4)$
can represent the presence/absence of a certain gene mutation. We used
two different sets of
values for $\beta$: the first set is
$\beta=(-3.45, 0.26, 0.1, 0.3, 0.26)^{\mathrm{T}}$, where $\beta_4=0.26$ represents a relatively common mutation;
the second set is $\beta=(-3.2, 0.26, 0.1, 0.3, 0.065)^{\mathrm{T}}$, where $\beta_4=0.065$ represents a very rare mutation.
In both sets, the true parameters are chosen so that the model
results in a population disease rate $p_D^t(1)\approx5\%$. The
simulation results are
presented in the upper half of Table~\ref{table:simu}.

%
\begin{table*}
\tabcolsep=0pt
\caption{Simulation results for the two experiments, each with two
different sets of parameter values, representing uncommon (upper-left)
and common (upper-right) gene mutation, and homogeneous (lower-left)
and diversified (lower-right) gene expression levels.
`true' is the true value of
$\beta$, `est' is the average of the estimated $\beta$, `sd' is the
sample standard deviation and `$\widehat{\operatorname{sd}}$' is the average of
the estimated standard deviation}\label{table:simu}
\begin{tabular*}{\textwidth}{@{\extracolsep{\fill}}ld{2.4}d{1.4}d{1.4}d{1.4}d{1.4} d{2.4}d{1.4}d{1.4}d{1.4}d{1.4}@{}}
\hline
&\multicolumn{1}{l}{$\beta_c$}
&\multicolumn{1}{l}{$\beta_1$}
&\multicolumn{1}{l}{$\beta_2$}
&\multicolumn{1}{l}{$\beta_3$}
&\multicolumn{1}{l}{$\beta_4$}
&\multicolumn{1}{l}{$\beta_c$}
&\multicolumn{1}{l}{$\beta_1$}
&\multicolumn{1}{l}{$\beta_2$}
&\multicolumn{1}{l}{$\beta_3$}
&\multicolumn{1}{l@{}}{$\beta_4$}\\
\hline
&\multicolumn{10}{c}{Experiment 1}\\
true & -3.2000 & 0.2600 & 0.1000 & 0.3000 & 0.0650 & -3.4500 & 0.2600 & 0.1000 & 0.3000 & 0.2600\\
est & -3.8925 & 0.2498 & 0.0995 & 0.3101 & 0.0649 & -3.9263 & 0.2618 & 0.0994 & 0.2998 & 0.2610\\
sd & 1.6390 & 0.3110 & 0.0359 & 0.1226 & 0.0111 & 1.3958 & 0.2196 & 0.0445 & 0.0783 & 0.0229\\
$\widehat{\operatorname{sd}}$ & 1.6285 & 0.3236 & 0.0364 & 0.1192 & 0.0116 & 1.2534 & 0.1956 & 0.0422 & 0.0723 & 0.0207\\[6pt]
 & \multicolumn{10}{c}{Experiment 2}\\
true & -3.2000 & 0.2600 & 0.1000 & 0.3000 & 0.3000 & -3.7300 & 0.2600 & 0.1000 & 0.3000 & 1.0000\\
est & -3.3128 & 0.2553 & 0.0993 & 0.3126 & 0.2999 & -3.7442 & 0.2589 & 0.0995 & 0.3053 & 0.9986\\
sd & 0.7815 & 0.1624 & 0.0352 & 0.0750 & 0.0101 & 0.2906 & 0.0685 & 0.0442 & 0.0405 & 0.0378\\
$\widehat{\operatorname{sd}}$ & 0.7969 & 0.1663 & 0.0358 & 0.0789 & 0.0101 & 0.2859 & 0.0676 & 0.0439 & 0.0402 & 0.0373\\
\hline
\end{tabular*}
\end{table*}

The second experiment differs from the first one in its assumption on
$q(g,\beta_4)$. Here, we model $q(g,\beta_4)$ with a Laplace
distribution with variance $\beta_4$. This kind of model is
typically used to model the gene expression level. To maintain an
approximate $5\%$ disease rate in the population, we used
$\beta=(-3.2, 0.26, 0.1, 0.3, 0.3)^{\mathrm{T}}$ and
$\beta=(-3.73, 0.26, 0.1, 0.3, 1)^{\mathrm{T}}$ as the true parameter
values. Again, in the first set, $\beta_4=0.3$
represents a small
variation in the population distribution for the gene expression
levels, resulting in a more homogeneous population in terms of this gene.
In the second set, $\beta_4=1$ represents a larger
variation, so the population is more diversified.
The simulation results are
presented in the lower half of Table~\ref{table:simu}.
In both experiments, 1000 simulations are implemented.

From Table~\ref{table:simu}, it is clear that the estimator for $\beta
$ is
consistent in all four situations and the estimated standard deviation
approximates the true one rather well. It is worth noting that the
first experiment is a repetition of the same setting as in Chatterjee
and Carroll (\citeyear{ChatterjeeCarroll2005}) and we observe very similar results.
Specifically, for $\beta_1,\beta_2, \beta_3, \beta_4$ in the
upper-left table,
their results for ``sd'' are 0.322, 0.037, 0.128, 0.0122, respectively,
and those in the upper-right table are 0.198, 0.043, 0.075 and 0.0273,
respectively. Although some numerical improvement can be observed in
certain parameters (for example $\beta_4$), it can be a result of
finite-sample performance and numerical issues. We conjecture that
the estimator in Chatterjee and Carroll (\citeyear{ChatterjeeCarroll2005}) is equivalent to the method proposed here, hence is
also efficient, although a rigorous proof is beyond the scope of this
paper. It is
also worth noting that the estimation of $\beta_c$ is more difficult
than the remaining components of $\beta$, in that the estimation
has large variability. This is especially prominent
in the discrete model setting for $q(g)$. Indeed, the estimation
result for $\beta_c$ has not been reported elsewhere and, without the
gene-environment independence, $\beta_c$ is known to be
unidentifiable (Prentice and Pyke (\citeyear{PrenticePyke1979})). This provides an intuitive
explanation for the performance of $\hat\beta_c$ we observe. The set of
estimating equations is solved using a standard Newton--Raphson algorithm.

\section{Conclusion}\label{sec:discuss}
Semiparametric modeling and estimation to study the occurrence of a
disease in relation to gene and environment has attracted much interest
recently. However, despite the various estimators proposed in
the literature, very little is understood in terms of the efficiency of
the estimators. This is partly due to the fact that most estimators
are constructed in rather ingenious ways, instead of following the standard
lines of semiparametric theory. The other reason is that most such
problems are set in a case-control design, which violates the
i.i.d.~assumption for standard semiparametric theory.

Instead of rederiving the whole semiparametric
theory under non-i.i.d.~samples, we argue that case-control data can be
treated as if they were i.i.d.~data and the standard
semiparametric efficiency theory will still apply. The equivalence of
the first order asymptotic theory shown in this article is a new contribution.
The argument is
based on rather elementary statistics without involving advanced
knowledge or highly
specialized techniques.

The establishment of the equivalence of the semiparametric efficiency
between i.i.d.~data and case-control data allows us to carry out the
estimation using standard, well-established semiparametric
theory. However, these standard analyses are performed under a
hypothetical population of interest, hence the detailed derivation often
requires special treatment, something which has not previously appeared in
the literature.
Under the gene-environment independence assumption, we are
able to explicitly construct a novel semiparametric estimator and show its
efficiency. A~special feature of this estimator is that we never
attempted to estimate the infinite-dimensional nuisance parameter $\eta$
itself, neither did we posit a model, true or false, for
it. Rather, we
avoided its estimation and instead approximated quantities
that rely on it. On the one hand, this enables us to carry out the
estimation rather easily; on the other hand, some asymptotic
properties have to be rederived because any result that relies on the
convergence properties of the nuisance parameter itself can no longer
be used.

Finally, our simulation results support the theory we developed,
in both discrete
and continuous gene distribution cases.
Our simulation results in the
discrete gene model are very similar to those of Chatterjee and Carroll
(\citeyear{ChatterjeeCarroll2005}), which leads us to believe that their estimator is also efficient.
A demonstration of this aspect would be an interesting direction for
future work.
The programming of the method in Chatterjee and
Carroll may seem easier. However, if the two methods are indeed
equivalent, then the projection step in the current method should be
equivalent to the profiling step in Chatterjee and Carroll, hence the
computational effort and intensity should be equivalent.
Although we did not further expand our estimator to stratified
case-control data,
the method is clearly applicable there as well.

\begin{appendix}
\section*{Appendix}\label{Appendix}
\setcounter{equation}{0}

\subsection*{The derivation of $S_{\mathrm{eff}}$}

We will use $S_{\mathrm{eff}}$ to construct our estimating
equation. We calculate $S_{\mathrm{eff}}$ by projecting the score
functions with respect to the parameters of interest $\beta_c,
\beta_1, \beta_2, \beta_3, \beta_4$ onto the orthogonal
complement of the nuisance tangent space.
We first derive the score functions $S_\beta\equiv\partial\log
p(g,e,d;\beta,\eta)/\partial\beta$. Straightforward calculation
shows that
the score function
$S_\beta=(S_1^{\mathrm{T}},S_2^{\mathrm{T}})^{\mathrm{T}}$, where
\begin{eqnarray*}
S_1^{\mathrm{T}}&=&(m'_{\beta_c}\ m'_{\beta_1}\ m'_{\beta_2}\ m'_{\beta_3})
\biggl(d-1+\frac{1}{1+\mathrm{e}^m}\biggr)
-E\biggl\{(m'_{\beta_c}\ m'_{\beta_1}\ m'_{\beta_2}\ m'_{\beta_3})
\biggl(d-1+\frac{1}{1+\mathrm{e}^m}\biggr) \Big|d\biggr\},\\
S_2&=&\frac{q'_{\beta_4}(g,{\beta_4})}{q(g,{\beta_4})}
-E\biggl\{\frac{q'_{\beta_4}(g,{\beta_4})}{q(g,{\beta_4})} \Big|d\biggr\}.
\end{eqnarray*}
Here, $m'_*$ and $q'_*$ represent partial derivatives with respect to
$*$.
Note that, in general, $S_\beta$ can be written as $S_\beta=S-E(S|d)$.

We next derive the nuisance tangent space $\Lambda$ and its orthogonal
complement $\Lambda^\perp$. Inserting the form of
$p_D^t(d;\beta,\eta)$ into (\ref{eq:pdf}), replacing $\eta(e)$ by
an arbitrary submodel $p^t_E(e;\gamma)$ and taking the derivative of
$\log
p(g,e,d;\beta,\gamma)$ with respect to $\gamma$, we obtain
$\partial\log p(g,e,d;\beta,\gamma)/\partial\gamma
=\partial\log p_E^t(e;\gamma)/\partial\gamma-
E\{\partial\log p_E^t(e;\gamma)/\partial\gamma|d\}$.
Now, recognizing that $\partial\log p_E^t(e;\gamma)/\partial\gamma$
for an arbitrary submodel can
yield an arbitrary function of $e$ with mean zero calculated under the true
$\eta(e)$, we obtain the nuisance
tangent space:
\begin{eqnarray*}
\Lambda&=&[h(e)-E\{h(e)|d\}{}\dvt{} \forall h(e) \mbox{ such that } E^t(h)=0]
=[h(e)-E\{h(e)|d\}{}\dvt{} \forall h(e)],\\
\Lambda^\perp&=&[h(g,e,d){}\dvt{} E(h|e)=E\{E(h|d)|e\}].
\end{eqnarray*}
Here, $E^t$ stands for an expectation calculated with respect to the true
population distribution.
The second expression for $\Lambda$ is more convenient because it
allows $h(e)$ to be an arbitrary function of~$e$, hence this is
the form of $\Lambda$ that we will use.

Having obtained $S_\beta$ and the spaces $\Lambda$ and
$\Lambda^\perp$, we can proceed to derive
the efficient score function
$S_{\mathrm{eff}}\equiv\Pi(S_\beta|\Lambda^\perp)$. If we let
$\Pi(S_\beta|\Lambda)=f(e)-E(f|d) $,
then
$S_{\mathrm{eff}}=S_\beta-f(e)+E(f|d)=S-E(S|d)-f(e)+E(f|d)$.

We now modify the expression of $S_{\mathrm{eff}}$ to facilitate
its actual
computation. Letting
$a(d)=E(f|d)-E(S|d)$, we can thus write
$S_{\mathrm{eff}}=S-f+a(d)$.
Note that $S$ does not depend on
$\eta$ and $a(d)$ is either $a(1)$ or $a(0)$.
In addition, we have
$E(S_{\mathrm{eff}}|e)=E\{E(S_{\mathrm{eff}}|d)|e\}$.
This is equivalent to
\begin{eqnarray*}
E(S_\beta|e)-f(e)+E\{E(f|d)|e\}
=E[E\{S-E(S|d)|d\}-E\{f-E(f|d)|d\}|e]=0,
\end{eqnarray*}
which, in turn, is equivalent to
\begin{eqnarray*}
E(S|e)&=&f+E\{E(S|d)|e\}-E\{E(f|d)|e\}=f-E\{a(d)|e\}\\
&=&f-\frac{\sum_d \int a(d)N_dq(g,{\beta_4})H(d,g,e)\,\mathrm{d}\mu(g)/p_D^t(d)}
{\sum_d\int N_dq(g,{\beta_4})H(d,g,e) \,\mathrm{d}\mu(g)/p_D^t(d)}.
\end{eqnarray*}
Let
\begin{eqnarray*}
v(e,d)=N_d\int
q(g,{\beta_4})H(d,g,e)\,\mathrm{d}\mu(g)/p_D^t(d)=p_{E,D}(e,d)N\eta^{-1}(e)
\end{eqnarray*}
and
%
\begin{eqnarray}
w(e,d)=v(e,d)/\{v(e,0)+v(e,1)\}.\label{eq:w}
\end{eqnarray}
We have
\begin{eqnarray*}
E(S|e)&=&f-a(0)v(e,0)/\{v(e,0)+v(e,1)\}-a(1)v(e,1)/\{v(e,0)+v(e,1)\}\\
&=&f-a(0)w(e,0)-a(1)w(e,1)
\end{eqnarray*}
or
$f=E(S|e)+a(0)w(e,0)+a(1)w(e,1)$.
Consequently,
\begin{eqnarray*}
S_{\mathrm{eff}}&=&S-E(S|e)-a(0)w(e,0)-a(1)w(e,1)+a(d)\\
&=&S-E(S|e)+(-1)^d\{a(0)-a(1)\}w(e,1-d).
\end{eqnarray*}

\subsection*{Proof of Theorem~\protect\ref{th:thm}}

To simplify notation, we denote $\alpha=p_D^t(0)/p_D^t(1)$,
$\hat\alpha=\hat p_D^t(0)/\hat p_D^t(1)$,
$\delta(\alpha)=a\{0;p_D^t(d)\}-a\{1;p_D^t(d)\}$,
$\delta(\hat\alpha)=a\{0;\hat p_D^t(d)\}-a\{1;\hat p_D^t(d)\}$ and
$\hat\delta(\hat\alpha)=\hat a\{0;\hat p_D^t(d)\}-\hat a\{1;\hat
p_D^t(d)\}$.

Suppose we randomly partition the data into two groups: group 1 has
$m$ observations and group 2 has $n$ observations.
Here, $m=N^{0.9}$, $n=N-m$.\vspace*{1pt} We use the first group to obtain $\hat
\alpha$,
and $\hat\delta(\hat\alpha)$,
then use the second group to form the following estimating equation to
estimate $\beta$:
\begin{eqnarray*}
\sum_{i=1}^n S_{\mathrm{eff}}\{x_i;\hat\beta,\hat\alpha
,\hat\delta(\hat
\alpha)\}=0.
\end{eqnarray*}
We will first show that the resulting estimator satisfies
$n^{1/2}(\hat\beta-\beta_0)\to N(0, V)$ in distribution when $N\to
\infty$.

The proof splits into several steps: First,
obviously, $\hat\alpha-\alpha=\mathrm{O}_p(m^{-1/2})$ and
$\hat\delta(\hat\alpha)-\delta(\hat\alpha)=\mathrm{O}_p(m^{-1/2})$, as
long as
a root-$N$-consistent $\tilde\beta$ is inserted in the calculation of these
quantities. A~standard expansion yields
\begin{eqnarray*}
0&=&\sum_{i=1}^n S_{\mathrm{eff}}\{x_i;\hat\beta,\hat\alpha,\hat\delta(\hat\alpha)\}\\
&=&\sum_{i=1}^n S_{\mathrm{eff}}\{x_i;\beta_0,\hat\alpha,\hat\delta(\hat\alpha)\}
+\sum_{i=1}^n\frac{\partial}{\partial\beta^{\mathrm{T}}}S_{\mathrm{eff}}\{x_i;\beta^*,\hat\alpha
,\hat\delta(\hat\alpha)\}(\hat\beta-\beta_0)
\\
&=&\sum_{i=1}^n S_{\mathrm{eff}}\{x_i;\beta_0,\hat\alpha,\hat\delta(\hat\alpha)\}
+n\biggl\{E\biggl(\frac{\partial S_{\mathrm{eff}}}{\partial\beta^{\mathrm{T}}}\biggr)+\mathrm{o}_p(1)\biggr\}(\hat\beta-\beta_0),
\end{eqnarray*}
which can be rewritten as
\begin{eqnarray*}
&&\biggl\{E\biggl(\frac{\partial S_{\mathrm{eff}} }{\partial\beta^{\mathrm{T}}}\biggr)+\mathrm{o}_p(1)\biggr\}
n^{1/2}(\hat\beta-\beta_0)
\\
&&\quad=-n^{-1/2}\sum_{i=1}^n
S_{\mathrm{eff}}\{x_i;\beta_0,\hat\alpha,\hat\delta(\hat\alpha)\}
\\
&&\quad=-n^{-1/2}\sum_{i=1}^n[
S_{\mathrm{eff}}\{x_i;\beta_0,\hat\alpha,\delta(\hat\alpha)\}
+(-1)^{d_i}\{\hat\delta(\hat\alpha)-\delta(\hat\alpha)\}
w(e_i,1-d_i,\hat\alpha)].
\end{eqnarray*}
The last equality\vspace*{1pt} uses the form of $S_{\mathrm{eff}}$ in (\ref
{eq:effscore})
and the fact that $S$, $E(S|e)$ and $w$ do not depend on $\alpha$.
Because
$\hat\delta(\hat\alpha)-\delta(\hat\alpha)=\mathrm{O}_p(m^{-1/2})=\mathrm{o}_p(1)$ and
\begin{eqnarray*}
E\{(-1)^{d_i}w(e_i,1-d_i,\hat\alpha)\}=\int\sum_{d=0,1}\frac
{(-1)^dp_{E,D}(e,1-d;\hat\alpha)\eta^{-1}(e)}{v(e,0;\hat\alpha
)+v(e,1;\hat\alpha)}p_{E,D}(e,d;\hat\alpha)\,\mathrm{d}\mu(e)=0,
\end{eqnarray*}
we actually have\vspace*{2pt}
\begin{eqnarray*}
&&\biggl\{E\biggl(\frac{\partial S_{\mathrm{eff}}}{\partial\beta^{\mathrm{T}}}\biggr)+\mathrm{o}_p(1)\biggr\}
n^{1/2}(\hat\beta-\beta_0)
\\
&&\quad=-n^{-1/2}\sum_{i=1}^n
S_{\mathrm{eff}}\{x_i;\beta_0,\hat\alpha,\delta(\hat\alpha)\}+\mathrm{o}_p(1)
\\
&&\quad=-n^{-1/2}\sum_{i=1}^n\biggl\{
S_{\mathrm{eff}}(x_i)+\frac{\partial
S_{\mathrm{eff}}(x_i;\beta_0,\alpha)}{\partial\alpha}(\hat
\alpha
-\alpha)
+\frac{\partial^2
S_{\mathrm{eff}}(x_i;\beta_0,\alpha^*)}{\partial\alpha^2}(\hat\alpha-\alpha)^2\biggr\}
+\mathrm{o}_p(1).
\end{eqnarray*}
In addition, $(\hat\alpha-\alpha)^2=\mathrm{O}_p(m^{-1})=\mathrm{o}_p(n^{-1/2})$, so\vspace*{2pt}
\begin{eqnarray*}
\biggl\{E\biggl(\frac{\partial S_{\mathrm{eff}} }{\partial\beta^{\mathrm{T}}}\biggr)+\mathrm{o}_p(1)\biggr\}
n^{1/2}(\hat\beta-\beta_0)
=-n^{-1/2}\sum_{i=1}^n\biggl\{
S_{\mathrm{eff}}(x_i)+\frac{\partial
S_{\mathrm{eff}}(x_i)}{\partial\alpha}(\hat\alpha-\alpha)\biggr\}
+\mathrm{o}_p(1).
\end{eqnarray*}

We now proceed to examine $\frac{\partial
S_{\mathrm{eff}}(x_i)}{\partial\alpha}$ by examining each term in
(\ref{eq:effscore}). $S$ is free of $\alpha$.
As a function of $\alpha$, we already have\vspace*{2pt}
\begin{eqnarray*}
b_5(e;\alpha)&\equiv& E(S|e;\alpha)
=\frac{\sum_d \int SN_dq(g,{\beta_4})H(d,g,e)
\,\mathrm{d}\mu(g)/p_D^t(d)}
{\sum_d\int N_dq(g,{\beta_4})H(d,g,e)\,\mathrm{d}\mu(g)/p_D^t(d)}\\
&=&
\frac{\int SN_0qH_0\,\mathrm{d}\mu(g)+\alpha\int SN_1qH_1\,\mathrm{d}\mu(g)}
{\int N_0qH_0\,\mathrm{d}\mu(g)+\alpha\int N_1qH_1\,\mathrm{d}\mu(g)}=
\frac{u_2(e,0)+\alpha u_2(e,1)}{u_1(e,0)+\alpha u_1(e,1)},
\end{eqnarray*}
where we define $u_1(e,d)=\int N_dq(g,{\beta_4})H(d,g,e)\,\mathrm{d}\mu(g)$ and
$u_2(e,d)=\int SN_dq(g,{\beta_4})H(d,\break g,e)\,\mathrm{d}\mu(g)$. Using this notation,\vspace*{2pt}
\begin{eqnarray*}
\frac{\partial b_5}{\partial\alpha}
&=&\frac{u_2(e,1)u_1(e,0)-u_2(e,0)u_1(e,1)}
{\{u_1(e,0)+\alpha u_1(e,1)\}^2},\\
w(e,0)&=&\frac{u_1(e,0)}{u_1(e,0)+\alpha u_1(e,1)},\\
w(e,1)&=&\frac{\alpha u_1(e,1)}{u_1(e,0)+\alpha u_1(e,1)}.
\end{eqnarray*}
Similarly to the calculation of $b_0, b_1$, we also have that for any
function $u$,\vspace*{2pt}
\begin{eqnarray*}
E(u|d;\alpha)
&=&\int\frac{p_E(e)\int
uN_dq(g)H(d,g,e)\,\mathrm{d}\mu(g)/p_D^t(d)}{\sum_d\int
N_dq(g)H(d,g,e) \,\mathrm{d}\mu(g)/p_D^t(d)}\,\mathrm{d}\mu(e)
\\
&&{}\Big/
\int\frac{p_E(e)\int
N_dq(g)H(d,g,e)\,\mathrm{d}\mu(g)/p_D^t(d)}{\sum_d\int
N_dq(g)H(d,g,e) \,\mathrm{d}\mu(g)/p_D^t(d)}\,\mathrm{d}\mu(e)\\
&=&\int\frac{p_E(e)\int
uN_dq(g)H(d,g,e)\,\mathrm{d}\mu(g)/p_D^t(d)}{\sum_d
u_1(e,d)/p_D^t(d)}\,\mathrm{d}\mu(e)
\\
&&{}\Big/
\int\frac{p_E(e)u_1(e,d)/p_D^t(d)}{\sum_du_1(e,d)/p_D^t(d)}\,\mathrm{d}\mu(e),
\end{eqnarray*}
thus
\begin{eqnarray*}
E(u|0;\alpha)=\int\frac{p_E(e)\int
uN_0q(g)H(0,g,e)\,\mathrm{d}\mu(g)}{u_1(e,0)+
u_1(e,1)\alpha}\,\mathrm{d}\mu(e) \Big/
\int\frac{p_E(e)u_1(e,0)}{u_1(e,0)+u_1(e,1)\alpha}\,\mathrm{d}\mu(e),\\
E(u|1;\alpha)=\int\frac{p_E(e)\int
uN_1q(g)H(1,g,e)\,\mathrm{d}\mu(g)}{u_1(e,0)+
u_1(e,1)\alpha}\,\mathrm{d}\mu(e) \Big/
\int\frac{p_E(e)u_1(e,1)}{u_1(e,0)+u_1(e,1)\alpha}\,\mathrm{d}\mu(e).
\end{eqnarray*}
These relations lead to
\begin{eqnarray*}
b_0&=&\int
\frac{p_E(e)u_2(e,0)}{u_1(e,0)+u_1(e,1)\alpha}\,\mathrm{d}\mu(e)
\Big/
\int
\frac{p_E(e)u_1(e,0)}{u_1(e,0)+u_1(e,1)\alpha}\,\mathrm{d}\mu(e),\\
b_1&=&\int
\frac{p_E(e)u_2(e,1)}{u_1(e,0)+u_1(e,1)\alpha}\,\mathrm{d}\mu(e)
\Big/
\int\frac{p_E(e)u_1(e,1)}{u_1(e,0)+u_1(e,1)\alpha}\,\mathrm{d}\mu(e),\\
b_2&\equiv&E\{E(S|e)|0;\alpha\}\\
&=&\int\frac{p_E(e)\int
E(S|e)N_0q(g)H(0,g,e)\,\mathrm{d}\mu(g)}{u_1(e,0)+
u_1(e,1)\alpha}\,\mathrm{d}\mu(e) \Big/
\int\frac{p_E(e)u_1(e,0)}
{u_1(e,0)+u_1(e,1)\alpha}\,\mathrm{d}\mu(e)\\
&=&\int\frac{p_E(e)E(S|e)u_1(e,0)}{u_1(e,0)+
u_1(e,1)\alpha}\,\mathrm{d}\mu(e) \Big/
\int\frac{p_E(e)u_1(e,0)}
{u_1(e,0)+u_1(e,1)\alpha}\,\mathrm{d}\mu(e)\\
&=&\int\frac{p_E(e)u_1(e,0)\{u_2(e,0)+u_2(e,1)\alpha\}}{\{u_1(e,0)+
u_1(e,1)\alpha\}^2}\,\mathrm{d}\mu(e) \Big/
\int\frac{p_E(e)u_1(e,0)}{u_1(e,0)+u_1(e,1)\alpha}\,\mathrm{d}\mu(e),\\
b_3&\equiv&E\{w(e,0)|D=0\}\\
&=&\int\frac{p_E(e)\int
w(e,0)N_0q(g)H(0,g,e)\,\mathrm{d}\mu(g)}{u_1(e,0)+
u_1(e,1)\alpha}\,\mathrm{d}\mu(e) \Big/
\int\frac{p_E(e)u_1(e,0)}{u_1(e,0)+u_1(e,1)\alpha}\,\mathrm{d}\mu(e)\\
&=&\int\frac{p_E(e)u_1^2(e,0)}{\{u_1(e,0)+
u_1(e,1)\alpha\}^2}\,\mathrm{d}\mu(e) \Big/
\int\frac{p_E(e)u_1(e,0)}{u_1(e,0)+u_1(e,1)\alpha}\,\mathrm{d}\mu(e).
\end{eqnarray*}
Consequently, we obtain
\begin{eqnarray*}
S_{\mathrm{eff}}(0)&=&S-b_5(e)+\biggl\{b_1-b_0+\frac{b_2-b_0b_3-b_1(1-b_3)}{1-b_3}\biggr\}
\frac{\alpha u_1(e,1)}{u_1(e,0)+\alpha u_1(e,1)}\\
&=&S-b_5(e)+\biggl(\frac{b_2-b_0}{1-b_3}\biggr)\frac{\alpha
u_1(e,1)}{u_1(e,0)+\alpha u_1(e,1)},\\
S_{\mathrm{eff}}(1)&=&S-b_5(e)-\biggl(\frac{b_2-b_0}{1-b_3}\biggr)\frac
{u_1(e,0)}{u_1(e,0)+\alpha u_1(e,1)},\\
\frac{\partial S_{\mathrm{eff}}(0)}{\partial\alpha}&=&-b_5(e)'
+\biggl(\frac{b_2-b_0}{1-b_3}\biggr)'\frac{\alpha u_1(e,1)}{u_1(e,0)+\alpha u_1(e,1)}
+\biggl(\frac{b_2-b_0}{1-b_3}\biggr)\frac{u_1(e,0)
u_1(e,1)}{\{u_1(e,0)+\alpha u_1(e,1)\}^2},\\
\frac{\partial S_{\mathrm{eff}}(1)}{\partial\alpha}&=&-b_5(e)'
-\biggl(\frac {b_2-b_0}{1-b_3}\biggr)'\frac{u_1(e,0)}{u_1(e,0)+\alpha u_1(e,1)}
+\biggl(\frac{b_2-b_0}{1-b_3}\biggr)\frac{u_1(e,0)
u_1(e,1)}{\{u_1(e,0)+\alpha u_1(e,1)\}^2}.
\end{eqnarray*}
Since $S$ does not contain $\alpha$,
$\frac{\partial S_{\mathrm{eff}}}{\partial\alpha}$ is a
function of $(e,d)$
only. Because $p_{E,D}(e,d)=\eta(e)u_1(e,d)/\{Np_D^t(d)\}$, we have
$p_{E,D}(e,0)=(1+\alpha)\eta(e)u_1(e,0)/(N\alpha)$,
$p_{E,D}(e,1)=(1+\alpha)\eta(e)u_1(e,1)/N$ and
$p_E(e)=(1+\alpha)\eta(e)\{u_1(e,0)+\alpha u_1(e,1)\}/(N\alpha)$.
Combining these results, we have
\begin{eqnarray*}
E\biggl(\frac{\partial S_{\mathrm{eff}}}{\partial\alpha}\biggr)
&=&E\biggl[-b_5'(e)+\biggl(\frac{b_2-b_0}{1-b_3}\biggr)\frac{u_1(e,0)
u_1(e,1)}{\{u_1(e,0)+\alpha u_1(e,1)\}^2}\biggr]\\
&=&E\biggl[\frac{-u_2(e,1)u_1(e,0)+u_2(e,0)u_1(e,1)}
{\{u_1(e,0)+\alpha u_1(e,1)\}^2}\biggr]
\\
&&{}+\biggl(\frac{b_2-b_0}{1-b_3}\biggr)E\biggl[\frac{u_1(e,0)
u_1(e,1)}{\{u_1(e,0)+\alpha u_1(e,1)\}^2}\biggr].
\end{eqnarray*}
Plugging in the expressions for $b_0, b_2, b_3$, we obtain
\begin{eqnarray*}
\frac{b_2-b_0}{1-b_3}
&=&\biggl[\int\frac{p_E(e)u_1(e,0)\{u_2(e,0)+u_2(e,1)\alpha\}}{\{u_1(e,0)+
u_1(e,1)\alpha\}^2}\,\mathrm{d}\mu(e)
-\int\frac{p_E(e)u_2(e,0)}{u_1(e,0)+u_1(e,1)\alpha}\,\mathrm{d}\mu(e)\biggr]\\
&&{}\Big/
\biggl[{\int\frac{p_E(e)u_1(e,0)}
{u_1(e,0)+u_1(e,1)\alpha}\,\mathrm{d}\mu(e)-\int\frac{p_E(e)u_1^2(e,0)}{\{u_1(e,0)+
u_1(e,1)\alpha\}^2}\,\mathrm{d}\mu(e)}\biggr]\\
&=&\int\frac{\alpha p_E(e)\{u_1(e,0)u_2(e,1)-u_1(e,1)u_2(e,0)\}}{\{u_1(e,0)+
u_1(e,1)\alpha\}^2}\,\mathrm{d}\mu(e)
\\
&&{}\Big/\Biggl[\int\frac{\alpha p_E(e)u_1(e,0)u_1(e,1)}{\{u_1(e,0)+u_1(e,1)\alpha\}^2}\,\mathrm{d}\mu(e)\Biggr]
\\
&=&
E\biggl[\frac{\{u_1(e,0)u_2(e,1)-u_1(e,1)u_2(e,0)\}}{\{u_1(e,0)+u_1(e,1)\alpha\}^2}\biggr]\Big/
E\biggl[\frac{u_1(e,0)u_1(e,1)}{\{u_1(e,0)+u_1(e,1)\alpha\}^2}\biggr],
\end{eqnarray*}
thus, we have $E(\partial S_{\mathrm{eff}}/\partial\alpha)=0$.

The fact that $E(\partial S_{\mathrm{eff}}/\partial\alpha)=0$, in
combination with
$\hat\alpha-\alpha=\mathrm{o}_p(1)$, yields
\begin{eqnarray*}
\biggl\{E\biggl(\frac{\partial S_{\mathrm{eff}} }{\partial\beta^{\mathrm{T}}}\biggr)+\mathrm{o}_p(1)\biggr\}
n^{1/2}(\hat\beta-\beta_0)
=-n^{-1/2}\sum_{i=1}^nS_{\mathrm{eff}}(x_i)+\mathrm{o}_p(1).
\end{eqnarray*}
Thus, we indeed have $n^{1/2}(\hat\beta-\beta_0)\sim N(0,V)$ asymptotically.

In fact, the classical $N^{1/2}(\hat\beta-\beta_0)\sim N(0,V)$ also holds.
This is because
\begin{eqnarray*}
N^{1/2}(\hat\beta-\beta_0)-n^{1/2}(\hat\beta-\beta_0)
=\frac{m}{N^{1/2}+n^{1/2}}(\hat\beta-\beta_0)
\to\frac{N^{0.9}}{N}n^{1/2}(\hat\beta-\beta_0)\to0
\end{eqnarray*}
when $N\to\infty$. Thus, our estimator is semiparametric efficient.
Because of the equivalence result developed in Section
\ref{sec:casecontrol}, the estimator is also semiparametric efficient
for case-control data.
We split the data set into two groups with sizes $m$ and $n$ for
simplicity of the asymptotic analysis. In reality, one can certainly
use the whole data set in each stage of the estimation.

\section*{Acknowledgments}
This work was supported by NSF Grant DMS-0906341.

\end{appendix}

\printhistory

\end{document}